\documentclass[11pt,twoside]{article}

\usepackage{epsfig,amsfonts,color}
\usepackage{amsmath}

\usepackage{amssymb, palatino, geometry,url}
\usepackage{algorithmic}
\usepackage[noresetcount,lined,boxed]{algorithm2e} 

\usepackage[colorlinks=true,linkcolor=blue,citecolor=blue,urlcolor=blue]{hyperref}
\usepackage{subcaption}
\usepackage{amsthm}
\newtheorem{thm}{Theorem}

\usepackage{fancyhdr}
\newcommand{\be}{\begin{equation}}
\newcommand{\ee}{\end{equation}}
\newcommand{\bea}{\begin{eqnarray}}
\newcommand{\eea}{\end{eqnarray}}

\newcommand{\bvec}{\left(\begin{array}{c}}
	\newcommand{\evec}{\end{array}\right)}
\newcommand{\bsub}{\begin{subequations}}
	\newcommand{\esub}{\end{subequations}}

\usepackage{lineno}

\title{A Mixed-Integer Conic Programming Formulation for\\ Computing the Flexibility Index under \\Multivariate Gaussian Uncertainty}
\author{Joshua L. Pulsipher and Victor M. Zavala\thanks{Corresponding Author: victor.zavala@wisc.edu}\\
	{\small Department of Chemical and Biological Engineering}\\
	{\small \;University of Wisconsin, 1415 Engineering Dr, Madison, WI 53706, USA}}
\date{}

\begin{document}

\maketitle

\begin{abstract}
We present a methodology for computing the flexibility index when uncertainty is characterized using  multivariate Gaussian random variables. Our approach computes the flexibility index by solving a mixed-integer conic program (MICP). This methodology directly characterizes ellipsoidal sets to capture correlations in contrast to previous methodologies that employ approximations. We also show that, under a Gaussian representation, the flexibility index can be used to obtain a lower bound for the so-called stochastic flexibility index (i.e., the probability of having feasible operation). Our results also show that the methodology can be generalized to capture different types of uncertainty sets. 
\end{abstract}


\section{Problem Definition and Setting}
 Flexibility seeks to quantify the ability of a physical system to counteract uncertainty (e.g. disturbances, parameters) in order to maintain feasible operation \cite{grossmann2014evolution}. Grossmann and co-workers have proposed diverse formulations and algorithms to assess system flexibility  \cite{grossmann2014evolution,swaney1985index1,swaney1985index2}. The so-called {\em flexibility test problem} seeks to find recourse variables $\mathbf{z} \in \mathbb{R}^{n_z}$ that counteract the uncertain parameters $\boldsymbol{\theta}\in T\subseteq \mathbb{R}^{n_\theta}$ in order to satisfy the system constraints $f_j( \mathbf{z}, \boldsymbol{\theta})\leq 0,\,j\in J$. The flexibility test problem is given by:
\begin{equation}
	\chi := \max_{\boldsymbol{\theta} \in T} \ \ \psi(\boldsymbol{\theta}).
	\label{eq:feasibility_definition}
\end{equation}
\noindent Here, $\psi(\boldsymbol{\theta})$ is a function that tests the feasibility of the system at a fixed value of the parameters $\boldsymbol{\theta}$. This function is evaluated by searching for a recourse $\mathbf{z}$ that minimizes the largest constraint violation:
\begin{equation}
	\psi(\boldsymbol{\theta}):= \min_{\mathbf{z}} \ \max_{j \in J} \ \ f_j( \mathbf{z}, \boldsymbol{\theta}).
	\label{eq:feasibility_function}
\end{equation}
\noindent  If $\psi(\boldsymbol{\theta}) \leq 0$ the system has feasible operation at $\boldsymbol{\theta}$ and thus $\chi \leq 0$ indicates that the system has feasible operation over the entire uncertainty set  $T$ (and is thus deemed to be flexible). 

Problem~\eqref{eq:feasibility_function} can be reformulated by using an upper-bounding variable $u\in\mathbb{R}$ as:
\begin{equation}
	\begin{aligned}
		&\psi(\boldsymbol{\theta}) = &\min_{\mathbf{z}, u} &&& \ \ u \\
		&&\text{s.t.} &&& \ \ f_j( \mathbf{z}, \boldsymbol{\theta}) \leq u, & j \in J.
	\end{aligned}
	\label{eq:mod_feasibility_function}
\end{equation}

The so-called {\em flexibility index problem} seeks to identify the largest uncertainty set $T(\delta)$ (where $\delta\in\mathbb{R}_+$ is a variable that scales the set) under which the system remains feasible. Formally, the flexibility index $F\in\mathbb{R}_+$ is defined as the solution of the problem:
\begin{equation}
	\begin{aligned}
		& F := & \max_{\delta\in\mathbb{R}_+} &&&\delta \\
		&&\text{s.t.} &&& \max_{\boldsymbol{\theta} \in T(\delta)} \ \ \psi(\boldsymbol{\theta}) \leq 0.
	\end{aligned}
	\label{eq:flexibility_index_definition}
\end{equation}

Parameterizing the uncertainty set $T(\delta)$ in terms of a single scalar variable $\delta$ is in general difficult. Most studies reported in the literature assume a {\em hyperbox} set of the form:
\begin{equation}
	T_{box}(\delta) = \{\boldsymbol{\theta} : \bar{\boldsymbol{\theta}} - \delta\Delta\boldsymbol{\theta}^- \leq \boldsymbol{\theta} \leq \bar{\boldsymbol{\theta}} + \delta\Delta\boldsymbol{\theta}^+\}
	\label{eq:hypercube}
\end{equation}
\noindent where $\bar{\boldsymbol{\theta}}$ is a nominal value and $\Delta\boldsymbol{\theta}^-,\Delta\boldsymbol{\theta}^+$ are maximum lower and upper deviations \cite{grossmann2014evolution}.

Problems \eqref{eq:feasibility_definition} and \eqref{eq:flexibility_index_definition} are conceptually attractive, but are computationally challenging due to their nested nature \cite{grossmann1983optimization}. Swaney and Grossmann proposed search procedures for the special case in which the solution of problems \eqref{eq:feasibility_definition} and \eqref{eq:flexibility_index_definition} correspond to vertices of $T_{box}$ and $T_{box}(\delta)$, as is the case when the system constraints $f_j(\cdot), j \in J,$ are convex \cite{swaney1985index1,swaney1985index2}. Vertex exploration suffers from computational limitations because an exponential number of vertices ($2^{n_\theta}$) must be evaluated \cite{grossmann2014evolution}. Grossmann and Floudas developed an alternative approach that uses an active-set approach to compute the flexibility index under hyperbox sets \cite{grossmann1987active}. For the case of linear constraints and hyperbox sets, this approach casts the flexibility index problem as a mixed-integer linear programming formulation.

Unfortunately, hyperbox representations of $T(\delta)$ do not adequately capture correlations in the parameters  \cite{rooney1999incorporating}.  To address this issue, Rooney and Biegler \cite{rooney1999incorporating} presented a methodology which considers ellipsoidal sets of the form:
\begin{equation}
	T_{ellip} = \{\boldsymbol{\theta} : (\boldsymbol{\theta} - \bar{\boldsymbol{\theta}})^T V_{\boldsymbol{\theta}}^{-1} (\boldsymbol{\theta} - \bar{\boldsymbol{\theta}}) \leq \chi^2_{n_p}(\alpha) \}
	\label{eq:hyperellipsoid}
\end{equation}
\noindent where $V_{\boldsymbol{\theta}}\in \mathbb{R}^{n_\theta\times n_{\theta}}$ is the covariance matrix of $\boldsymbol{\theta}$ (assumed to be symmetric positive definite), $n_\theta$ is the dimension of the uncertainty region, and $\mathcal{\chi}^2_{n_\theta}(\alpha)$ is the critical value of a $\chi$-squared distribution with $n_\theta$ degrees of freedom and at probability level $\alpha$. Ellipsoidal uncertainty sets can be constructed from the level sets of multivariate Gaussian random variables, and thus can be used to capture correlations. In the approach of Rooney and Biegler, the uncertainty set is discretized along the longest axis and the discrete points are used in combination with the flexibility test problem~\eqref{eq:MILP_flexibility_test} to assess the flexibility of the system. This approach is intuitive and practical, but only represents an approximation \cite{rooney1999incorporating}. 

To capture more general cheracterizations of uncertainty, Straub and Grossmann developed a fundamentally different approach to flexibility analysis. This was done in terms of what they called the {\em stochastic flexibility index} \cite{straub1990integrated}:
\begin{equation}
	SF := \int_{\boldsymbol{\theta}\in\Theta} p(\boldsymbol{\theta}) d\boldsymbol{\theta}
	\label{eq:SF_definition}
\end{equation}
\noindent where $\Theta:=\{\boldsymbol{\theta}\,|\,\psi(\boldsymbol{\theta}) \leq 0\}$ is the feasible set (projected onto the uncertainty space) and the random parameter $\boldsymbol{\theta}$ is characterized by the joint probability density function $p:\mathbb{R}^{n_\theta}\to \mathbb{R}$. Pistikopoulos and Mazzuchi proposed a similar definition in \cite{pistikopoulos1990novel}. The stochastic flexibility index is interpreted as the probability that the system remains feasible. More specifically, it represents the probability of finding a recourse variable $z$ such that that the system constraints $f_j(\mathbf{z},\boldsymbol{\theta})\leq 0,\, j \in J$ are satisfied. Consequently, we observe that this index shares common ground with joint chance constraints, in which one seeks to enforce the system constraints with a desired level of probability. In particular, it is not difficult to see that the stochastic flexibility index satisfies $SF=\mathbb{P}(\psi(\boldsymbol{\theta}) \leq 0)=F_{\psi(\boldsymbol{\theta})}(0)$, where $F_{\psi(\boldsymbol{\theta})}:\mathbb{R}\to[0,1]$ is the cumulative density function of the random variable $\psi(\boldsymbol{\theta})=\min_{z} \ \max_{j \in J} \ \ f_j( \mathbf{z}, \boldsymbol{\theta})$.

Index $SF$ can be obtained directly by integrating $p(\cdot)$ over the feasible region projected in the $\boldsymbol{\theta}$ space \cite{grossmann2014evolution}. This can be done using Monte Carlo (MC) sampling and requires checking for feasibility at every sampled realization. It is well-known that MC has asymptotic convergence properties but might require an extremely large number of samples to cover the uncertainty space \cite{shapiro2013sample, robert2013monte}. This is particularly important in the context of flexibility analysis because limiting behavior is often found in the boundaries of the feasible region. To circumvent these issues, Straub and Grossmann proposed a quadrature scheme \cite{straub1990integrated}, which explores the random space in a more systematic manner and thus require far fewer samples than MC. Significant advances in quadrature schemes have been achieved in recent years due to the development of sparse grid techniques \cite{ma1996generalized}. However, quadrature still suffers from severe computational limitations in high dimensions \cite{gerstner1998numerical, straub1990integrated}. 

This paper proposes an approach to compute the flexibility index when uncertainty can be represented using multivariate Gaussian random variables. We show that the problem can be cast as a mixed-integer conic program (MICP) that seeks to find the largest radius of an ellipsoidal uncertainty set under which the system maintains feasible operation. This approach is based on the key observation that the active-set approach of Grossmann and Floudas can be applied to general uncertainty sets that can be parameterized using a scalar variable $\delta$. The proposed approach offers the ability to capture correlations and to provide probabilistic estimates of feasible operation. In particular, we show that the flexibility index can be used to compute a lower bound for the stochastic flexibility index. 

The paper is structured as follows. In Section \ref{sec:active_set} we review fundamental properties of active-set formulations of the flexibility index to argue that this can be applied to general uncertainty sets. In Section \ref{sec:micp} we propose the mixed-integer conic formulation for ellipsoidal sets and establish its properties. Illustrative examples are provided in Section \ref{sec:examples}. The paper closes with concluding remarks and directions of future work. 

\section{Mixed-Integer Formulations of Flexibility Index}
\label{sec:active_set}

Grossmann and Floudas presented a mixed-integer programming approach for solving the flexibility test and index problems~\eqref{eq:feasibility_definition} and~\eqref{eq:flexibility_index_definition} \cite{grossmann1987active}. This approach replaces the inner minimization problem \eqref{eq:mod_feasibility_function} by its first-order Karush-Kuhn-Tucker (KKT) conditions and reformulates the complementarity conditions using binary variables:
\begin{nolinenumbers}
	\begin{subequations}
		\begin{equation}
			\sum_{j \in J} \lambda_j = 1
			\label{eq:KKT_1}
		\end{equation}
		\begin{equation}
			\sum_{j \in J} \lambda_j \frac{\partial f_j(\mathbf{z},\boldsymbol{\theta})}{\partial \mathbf{z}} = 0
			\label{eq:partial}
		\end{equation}
		\begin{equation}
			\lambda_j(f_j(\mathbf{z},\boldsymbol{\theta})- u) = 0 \ \ \ j \in J
			\label{eq:KKT_activator}
		\end{equation}
		\begin{equation}
			f_j(\mathbf{z},\boldsymbol{\theta})- u + s_j = 0 \ \ \ j \in J
			\label{eq:KKT_system}
		\end{equation}
		\begin{equation}
			\lambda_j \geq 0, \ s_j \geq 0 \ \ \ j \in J
		\end{equation}
		\label{eq:KKT_conditions}
	\end{subequations}
\end{nolinenumbers}
\noindent where $\lambda_j$ are the Lagrange multipliers of the system constraints and $s_j$ are the corresponding slack variables. When the system constraints are convex, the KKT conditions are necessary and sufficient. Otherwise, the first-order conditions can also be satisfied by saddle points and maxima. Note that, since this transformation is applied for a fixed value of $\boldsymbol{\theta}$, no restrictions on the shape of the uncertainty set are imposed in this derivation. Because the optimal solution of \eqref{eq:mod_feasibility_function} occurs at $\psi(\boldsymbol{\theta}) = u$, the feasibility test problem can be written as:
\begin{equation}
	\begin{aligned}
		&\chi = &\max_{\boldsymbol{\theta} \in T} &&& u \\
		&&\text{s.t.} &&& \sum_{j \in J} \lambda_j = 1 \\
		&&&&& \sum_{j \in J} \lambda_j \frac{\partial f_j(\mathbf{z},\boldsymbol{\theta})}{\partial \mathbf{z}} = 0 \\
		&&&&& \lambda_j(f_j(\mathbf{z},\boldsymbol{\theta})- u) = 0 && j \in J \\
		&&&&& f_j(\mathbf{z},\boldsymbol{\theta})- u + s_j = 0 && j \in J \\
		&&&&& \lambda_j \geq 0, \ s_j \geq 0 && j \in J.
	\end{aligned}
	\label{eq:KKT_feasibility_function}
\end{equation}
\noindent  Problem~\eqref{eq:KKT_feasibility_function} is equivalent to the feasibility constraint in the flexibility index problem \eqref{eq:flexibility_index_definition}. The complementarity conditions~\eqref{eq:KKT_activator} are nonlinear expressions that are difficult to handle computationally. To avoid this issue, we define binary variables $y_j \in \{0, 1\}$ and write the system of constraints: 
\begin{equation}
	\begin{aligned}
		s_j &\leq U(1 - y_j) & j \in J \\
		\lambda_j &\leq y_j & j \in J
	\label{eq:slack_binary_relation}
	\end{aligned}
\end{equation}
\noindent where $U$ is a suitable upper bound for the slack variables $s_j$.  

Previous work by Swaney and Grossmann has established conditions under which there exist exactly $n_z + 1$ active constraints~\cite{swaney1985index2}. This remarkable result is formalized in Theorem \ref{thm:active}.  

\begin{thm}
	\label{thm:active}
	If the set of gradients $\frac{\partial}{\partial \mathbf{z}} [f_j(\mathbf{z},\boldsymbol{\theta})], j \in J,$ are linearly independent, then there exists $n_z + 1$ active constraints $f_j(\boldsymbol{\theta}, \mathbf{z}) \leq 0$ at the solution of \eqref{eq:KKT_feasibility_function}.
\end{thm}

A proof of Theorem \ref{thm:active} is provided in Appendix \ref{a:active}. This is a summary of the proof provided in \cite{swaney1985index1} (which has been adapted to our context and introduced here for completeness).  Theorem \ref{thm:active} can be used to bound the binary variables as:
\begin{equation}
	\sum_{j \in J} y_j = n_z + 1.
	\label{eq:nz_constraint}
\end{equation}
This constraint can greatly facilitate the search of active constraints. 

We now let \eqref{eq:slack_binary_relation} and \eqref{eq:nz_constraint} replace the constraints in~\eqref{eq:KKT_activator} and~\eqref{eq:KKT_system} in \eqref{eq:KKT_feasibility_function} to produce the following mixed-integer formulation of the flexibility test problem:
\begin{equation}
	\begin{aligned}
		&\chi = &\max_{u, \boldsymbol{\theta}, \mathbf{z}, \lambda_j, s_j, y_j} &&& u \\
		&&\text{s.t.} &&& f_j(\mathbf{z},\boldsymbol{\theta})+ s_j = u && j \in J \\
		&&&&& \sum_{j \in J} \lambda_j = 1 \\
		&&&&& \sum_{j \in J} \lambda_j \frac{\partial f_j(\mathbf{z},\boldsymbol{\theta})}{\partial \mathbf{z}} = 0 \\
		&&&&& s_j \leq U(1 - y_j) && j \in J \\
		&&&&& \lambda_j \leq y_j && j \in J \\
		&&&&& \sum_{j \in J} y_j = n_z + 1  \\
		&&&&& \boldsymbol{\theta} \in T \\
		&&&&& \lambda_j, s_j \geq 0,y_j \in \{0, 1\} && j \in J.
	\end{aligned}
	\label{eq:MILP_flexibility_test}
\end{equation}
When the system constraints are linear, the feasibility test problem is a mixed-integer linear program (MILP). 

Swaney and Grossmann also proved that the flexibility index problem can be expressed as the minimum scaling value $\delta$ satisfying $\psi(\boldsymbol{\theta}) = 0$. Consequently, \eqref{eq:flexibility_index_definition} can be formulated as:
	\begin{equation}
		\begin{aligned}
			& F = & \min_{\delta \in \mathbb{R}_+, \ \boldsymbol{\theta} \in T(\delta)} &&&\delta \\
			&&\text{s.t.} &&& \psi(\boldsymbol{\theta}) = 0.
		\end{aligned}
		\label{eq:flexibility_index_mod}
	\end{equation}
This result is established in the following theorem and was originally proved in \cite{swaney1985index1}. A summary of the proof  is provided in Appendix \ref{a:flexibility_index}.

\begin{thm}
	\label{thm:flexibility_index}
	If $T(\delta)$ is a compact set and the system constraints $f_j(\mathbf{z},\boldsymbol{\theta}), j \in J,$ are continuous in $\mathbf{z}$ and $\boldsymbol{\theta}$, the flexibility index problem can be written as \eqref{eq:flexibility_index_mod}.
\end{thm}

This result implies that the constraints of \eqref{eq:MILP_flexibility_test} can be substituted into problem~\eqref{eq:flexibility_index_mod}. This gives the mixed-integer formulation: 
\begin{equation}
	\begin{aligned}
		&F = &\min_{\delta, \boldsymbol{\theta}, \mathbf{z}, \lambda_j, s_j, y_j} &&& \delta \\
		&&\text{s.t.} &&& f_j(\mathbf{z},\boldsymbol{\theta})+ s_j = 0 && j \in J \\
		&&&&& \sum_{j \in J} \lambda_j = 1 \\
		&&&&& \sum_{j \in J} \lambda_j \frac{\partial f_j(\mathbf{z},\boldsymbol{\theta})}{\partial \mathbf{z}} = 0 \\
		&&&&& s_j \leq U(1 - y_j) && j \in J \\
		&&&&& \lambda_j \leq y_j && j \in J \\
		&&&&& \sum_{j \in J} y_j = n_z + 1  \\
		&&&&& \boldsymbol{\theta} \in T(\delta) \\
		&&&&& \lambda_j, s_j \geq 0,y_j \in \{0, 1\} && j \in J.
	\end{aligned}
	\label{eq:MILP_flexibility_index}
\end{equation}
\noindent When the system constraints are linear, this problem is also a MILP. 

A key observation is that Theorems \ref{thm:active} and \ref{thm:flexibility_index} hold for any compact set $T(\delta)$. A hyperbox representation of $T$ and $T(\delta)$ (denoted $T_{box}$ and $T_{box}(\delta)$) is often used because this yields a MILP formulation (provided that the constraints are linear). As we have discussed, however, hyperbox representations cannot adequately capture correlation information.

\section{Mixed-Integer Conic Formulation for Multivariate Gaussian Uncertainty}\label{sec:micp}

We propose an approach to compute the flexibility index when uncertainty is represented as a multivariate Gaussian variable $\boldsymbol{\theta}\sim \mathcal{N}(\bar{\boldsymbol{\theta}}, V_{\boldsymbol{\theta}})$. Although ellipsoidal regions have been explored in the literature, no approach has been proposed to identify the largest ellipsoidal region $T_{ellip}(\delta)$ for which feasible operation can be achieved. We now proceed to show that this problem can be cast as a mixed-integer conic program and that the resulting flexibility index provides an alternative metric that captures parameter correlations. We will also show that this new index can be used to obtain a lower bound for the stochastic flexibility index.

We consider an ellipsoidal set of the form:
\begin{equation}
	T_{ellip}(\delta) = \{\boldsymbol{\theta} : (\boldsymbol{\theta} - \bar{\boldsymbol{\theta}})^T V_{\boldsymbol{\theta}}^{-1} (\boldsymbol{\theta} - \bar{\boldsymbol{\theta}}) \leq \delta \}
	\label{eq:hyperellipsoid2}
\end{equation}
\noindent where $\delta\in \mathbb{R}_+$ is the radius of the ellipsoid (this variable is used to scale the set). The flexibility index problem that we propose  thus seeks to find the largest radius $\delta$ for which feasible operation can be obtained. As we have seen, any compact set can be embedded in the flexibility index problem \eqref{eq:MILP_flexibility_index}. We thus embed \eqref{eq:hyperellipsoid2} into the flexibility index problem to obtain:
\begin{equation}
	\begin{aligned}
		&\delta^* = &\min_{\delta, \mathbf{z}, \boldsymbol{\theta}, \lambda_j, s_j, y_j} &&& \delta \\
		&&\text{s.t.} &&& f_j(\mathbf{z},\boldsymbol{\theta})+ s_j = 0 && j \in J \\
		&&&&& \sum_{j \in J} \lambda_j = 1 \\
		&&&&& \sum_{j \in J} \lambda_j \frac{\partial f_j(\mathbf{z},\boldsymbol{\theta})}{\partial \mathbf{z}} = 0 \\
		&&&&& s_j \leq U(1 - y_j) && j \in J \\
		&&&&& \lambda_j \leq y_j && j \in J \\
		&&&&& \sum_{j \in J} y_j = n_z + 1  \\
		&&&&& (\boldsymbol{\theta} - \bar{\boldsymbol{\theta}})^T V_{\boldsymbol{\theta}}^{-1} (\boldsymbol{\theta} - \bar{\boldsymbol{\theta}}) \leq \delta  \\
		&&&&& \lambda_j, s_j \geq 0 ; \ \ \ y_j \in \{0, 1\} && j \in J.
	\end{aligned}
	\label{eq:MIQCP_flexibility_index}
\end{equation}

Problem~\eqref{eq:MIQCP_flexibility_index} is a mixed-integer conic program that directly characterizes the ellipsoidal uncertainty region via a conic constraint, thus avoiding the sub-optimality problems associated with ellipsoidal approximations \cite{rooney1999incorporating}. The recent emergence of commercial MICP solvers can be exploited to make the solution of Problem~\eqref{eq:MIQCP_flexibility_index} more computationally attractive \cite{berthold2012extending}. 

Now we establish properties of the solution to Problem \eqref{eq:MIQCP_flexibility_index}. We denote the solution of the problem as $(\delta^*,\boldsymbol{\theta}^*)$. We recall the feasible set of the system is given by $\Theta=\{\boldsymbol{\theta}\,|\,\psi(\boldsymbol{\theta})\leq 0\}$ and note that its boundary is given by $\partial\Theta=\{\boldsymbol{\theta}\,|\,\psi(\boldsymbol{\theta})=0\}$. We also note that the boundary of the uncertainty set $T_{ellip}(\delta^*)$ is given by $\partial T_{ellip}(\delta^*)=\{\boldsymbol{\theta}\,|\,(\boldsymbol{\theta}-\bar{\boldsymbol{\theta}})^TV_{\boldsymbol{\theta}^{-1}}(\boldsymbol{\theta}-\bar{\boldsymbol{\theta}})= \delta^*\}$. We recall that Theorem \ref{thm:flexibility_index} proves that the critical parameter $\boldsymbol{\theta}^*$ lies in the boundary of the feasible set (i.e., $\boldsymbol{\theta}^*\in\partial \Theta$).  We now proceed to show that the critical point also lies at the boundary of the uncertainty set $T_{ellip}(\delta^*)$ and that the entire uncertainty set $T_{ellip}(\delta^*)$ lies inside the feasible region $\Theta$. 

\begin{thm}\label{thm:solution_properties} The solution pair $(\delta^*,\boldsymbol{\theta}^*)$ satisfies the following properties: i) The uncertainty set is contained in the feasible set ($T_{ellip}(\delta^*)\subseteq\Theta$),  ii) the critical parameter $\boldsymbol{\theta}^*$ lies on the boundary of the uncertainty set ($\boldsymbol{\theta}^*\in T_{ellip}(\delta^*)$), and iii) the critical parameter $\boldsymbol{\theta}^*$  lies in the intersection of the boundaries of the feasible and uncertainty sets ($\boldsymbol{\theta}^*\in \partial \Theta\cap \partial T_{ellip}(\delta^*)$).
\end{thm}

\begin{proof} To prove i) assume there exists $\tilde{\boldsymbol{\theta}}\in T_{ellip}(\delta^*)$ but $\tilde{\boldsymbol{\theta}}\notin \Theta$. Since $\tilde{\boldsymbol{\theta}}$ is infeasible, we have that $\psi(\tilde{\boldsymbol{\theta}})>0$ which implies that $\psi(\tilde{\boldsymbol{\theta}})>\psi(\boldsymbol{\theta}^*)$ because $\psi(\boldsymbol{\theta}^*)=0$. This is a contradiction because $\boldsymbol{\theta}^* \in  \textrm{argmax}_{\boldsymbol{\theta}\in T_{ellip}(\delta^*)}\psi(\boldsymbol{\theta})$ and thus $\psi(\boldsymbol{\theta}^*)=0$ is the maximum possible value of $\psi(\boldsymbol{\theta})$ in  $ T_{ellip}(\delta^*)$.  To prove ii) we  note that $\boldsymbol{\theta}^*\in \partial\Theta$ holds and, therefore, if $\boldsymbol{\theta}^*\notin \partial T_{ellip}(\delta^*)$ then there exists $\tilde{\boldsymbol{\theta}}\in T_{ellip}(\delta^*)$ with $\tilde{\boldsymbol{\theta}}\notin \Theta$. The result then follows from the argument used to prove i). The proof of iii) follows trivially from the observation that $\boldsymbol{\theta}^*$ lies on the boundary of both $\Theta$ and $T_{ellip}(\delta^*)$.
\end{proof}

Given the solution $\delta^*$ of the flexibility index problem \eqref{eq:MIQCP_flexibility_index}, we can also compute the following quantity (that we refer to as the {\em confidence level}):
\begin{equation}
	\alpha^* := F_{n_\theta}(\delta^*)
	\label{eq:alpha_definition}
\end{equation}
\noindent where $F_{n_\theta}(\delta^*)$ denotes the cumulative density of a $\chi$-squared distribution with $n_\theta$ degrees of freedom and at critical value $\delta^*$. The cumulative density can be computed from:
\begin{align}
F_{n_\theta}(\delta^*)=\frac{\gamma(\frac{n_\theta}{2}, \frac{\delta^*}{2})}{\Gamma(\frac{n_\theta}{2})}
\end{align}
where $\gamma(\cdot)$ and $\Gamma(\cdot)$ are the incomplete and complete gamma functions \cite{nist}. The reasoning behind the definition of the confidence level becomes clear in the following theorem. 

\begin{thm}\label{thm:main} The flexibility index $\delta^*$ and associated confidence level $\alpha^*$ satisfy the following properties: i) $\alpha^*=\mathbb{P}(\boldsymbol{\theta}\in T_{ellip}(\delta^*))$,  ii) $\alpha^*$ is a lower bound for the stochastic flexibility index (i.e., $\alpha^*\leq SF$), and iii) $\delta^*$ is the critical value (quantile) associated with confidence level $\alpha^*$. 
\end{thm}

\begin{proof}
Because the covariance matrix is symmetric positive definite, we can apply an eigenvalue decomposition of the form $V_{\boldsymbol{\theta}}^{-1}=Q\Lambda Q^{-1}$ to express the ellipsoidal constraint $(\boldsymbol{\theta} - \bar{\boldsymbol{\theta}})^T V_{\boldsymbol{\theta}}^{-1} (\boldsymbol{\theta} - \bar{\boldsymbol{\theta}}) \leq \delta$ as $\mathbf{x}^T\mathbf{x} \leq \delta$ or $\sum_{i=1}^{n_\theta}x_i^2 \leq \delta$. Here, 
$\mathbf{x}=\sqrt{\Lambda} Q^{-1}(\boldsymbol{\theta} - \bar{\boldsymbol{\theta}})$, $\Lambda\in\mathbb{R}^{n_\theta\times n_\theta}$ is a diagonal matrix containing the eigenvalues of the covariance matrix, and where $Q\in \mathbb{R}^{n_\theta\times n_\theta}$ is an orthogonal matrix that contains the corresponding eigenvectors (thus $Q^{-1} = Q^T$). We also have that  $x_i\sim\mathcal{N}(0,1)$ for all $i=1,...,n_\theta$ (this can be shown by using the property that $\mathbf{x}$ is a linear transformation of the Gaussian random variable $\boldsymbol{\theta}$ and thus is also Gaussian). Consequently, $\mathbf{x}^T\mathbf{x}=\sum_{i=1}^{n_\theta}x_i^2$ follows a $\chi$-squared distribution with $n_\theta$ degrees of freedom and $\mathbb{P}(\boldsymbol{\theta}\in T_{ellip}(\delta))=\mathbb{P}(\mathbf{x}^T\mathbf{x}\leq \delta)$.  The cumulative density of a $\chi$-squared random variable  at critical value $\delta$ is given by \eqref{eq:alpha_definition}. Consequently, we have that $\alpha^*=\mathbb{P}(\boldsymbol{\theta}\in T_{ellip}(\delta^*))$, establishing i).  

To prove ii), we recall  from Theorem \ref{thm:solution_properties} that $T_{ellip}(\delta^*)\subseteq\Theta$ and thus:
\begin{align*}
\alpha^*&=\mathbb{P}(\boldsymbol{\theta}\in T_{ellip}(\delta^*))\\
&=\int_{\boldsymbol{\theta}\in T(\delta^*)}p(\boldsymbol{\theta})d\boldsymbol{\theta}\nonumber\\
&\leq \int_{\boldsymbol{\theta}\in\Theta}p(\boldsymbol{\theta})d\boldsymbol{\theta}\nonumber\\
&=\mathbb{P}(\boldsymbol{\theta}\in \Theta)\\
&=SF, 
\end{align*}
Consequently, $\alpha^*$ provides a lower bound for the stochastic flexibility index. To prove iii), we note that inversion of the cumulative density function \eqref{eq:alpha_definition} is given by its quantile (critical value), which is given by $\delta^*$. 
\end{proof}

{\bf Remarks:} As noted in the proof of Theorem \ref{thm:main}, the conic constraint $(\boldsymbol{\theta} - \bar{\boldsymbol{\theta}})^T V_{\boldsymbol{\theta}}^{-1} (\boldsymbol{\theta} - \bar{\boldsymbol{\theta}}) \leq \delta$ can be expressed as a sum of squares constraint $\displaystyle\sum_{i=1}^{n_\theta}x_i^2 \leq \delta$ by using an eigenvalue decomposition of the covariance matrix. This can facilitate numerical implementation. We also note that the shapes of the feasible region $\Theta$ and of the uncertainty set $T_{ellip}(\delta^*)$ affect the gap between the confidence level $\alpha^*$ and the stochastic flexibility index $SF$.

\section{Illustrative Examples}\label{sec:examples}

We illustrate the concepts by using a couple of examples. All of the formulations are implemented in JuMP 0.18.2 \cite{DunningHuchetteLubin2017} and are solved using Gurobi 7.5.2 on an Intel\textregistered \, Core\texttrademark \, i7-7500U machine running at 2.90 GHz with 4 hardware threads and 16 GB of RAM running Windows 10 Home.

\subsection{Simple System} 
Consider a system whose feasible region is described by four constraints. The system is subjected to two Gaussian random variables $\boldsymbol{\theta}$ and has no recourse variables $\mathbf{z}$:
\begin{equation}
	\begin{aligned}
		f_1 &= \theta_1 + \theta_2 - 14 \leq 0 \\
		f_2 &= \theta_1 - 2 \theta_2 - 2 \leq 0 \\
		f_3 &= -\theta_1 \leq 0 \\
		f_4 &= -\theta_2 \leq 0 
	\end{aligned}
	\label{eq:2d_example_constraints}
\end{equation}
\noindent We let the mean $\bar{\boldsymbol{\theta}}$ and the covariance matrix $V_{\boldsymbol{\theta}}$ be given by 
\begin{equation}
	\bar{\boldsymbol{\theta}} =
	\begin{bmatrix}
		4 \\ 5
		\end{bmatrix}
		\ \ \ \ \ \
		V_{\boldsymbol{\theta}} = 
		\begin{bmatrix}
		2 & \beta \\
		\beta & 3
	\end{bmatrix}
	\label{eq:2d_example_conditions}
\end{equation}
\noindent where $\beta=\textrm{Cov}(\theta_1,\theta_2)$ is the covariance between $\theta_1$ and $\theta_2$. Three cases are considered where the value of $\beta$ is taken to be -1, 0, and 1, respectively. In each case, the confidence level $\alpha^*$ is verified via MC sampling  by drawing 10,000 realizations from $\mathcal{N}(\bar{\boldsymbol{\theta}}, V_{\boldsymbol{\theta}})$ and by counting instances that lie within $T_{ellip}(\delta^*)$.  

The stochastic flexibility index is determined also via MC sampling by determining the feasibility of 100,000 realizations drawn from $\mathcal{N}(\bar{\boldsymbol{\theta}}, V_{\boldsymbol{\theta}})$ using Problem \eqref{eq:mod_feasibility_function}. This $SF$-MC sampling requires $421$-$448$ seconds to accurately estimate $SF$ in each case, while the proposed MICP formulation only requires $0.14$-$0.18$ seconds to compute $\alpha^*$.

Table~\ref{tab:2d_results} presents the flexibility index $\delta^*$ and confidence level $\alpha^*$ for each case. We can see that the $\alpha^*$-MC estimates are in agreement with the MICP results. We also confirm that the confidence level provides a lower bound for the stochastic flexibility index.

\begin{table}[ht]
	\caption{Results for system governed by~\eqref{eq:2d_example_constraints} and~\eqref{eq:2d_example_conditions}.}
	\begin{center}
		\begin{tabular}{l|ccccc}
			& $\beta$ & $\delta^*$ & $\alpha^*$-MICP (\%) & $\alpha^*$-MC (\%) & $SF$-MC (\%) \\ \hline
			Case 1 & -1 & 3.56 & 83.1 & 83.3 & 96.6\\
			Case 2 & 0 & 4.57 & 89.8 & 90.0 & 96.9\\
			Case 3 & 1 & 3.57 & 83.2 & 83.1 & 96.3 
		\end{tabular}
	\end{center}
	\label{tab:2d_results}
\end{table}

\begin{figure}[ht]
	\centering
	\begin{subfigure}[t]{.3\textwidth}
		\centering
		\includegraphics[height=2.3in]{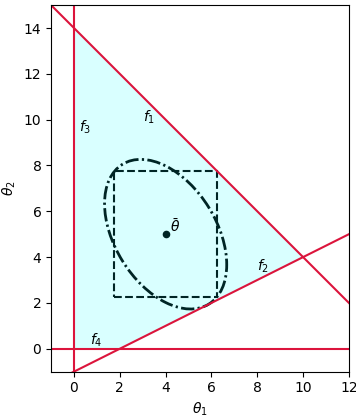}
		\caption{Case 1: $\textrm{Cov}(\theta_1,\theta_2) = -1$}
	\end{subfigure}
	\begin{subfigure}[t]{.3\textwidth}
		\centering
		\includegraphics[height=2.3in]{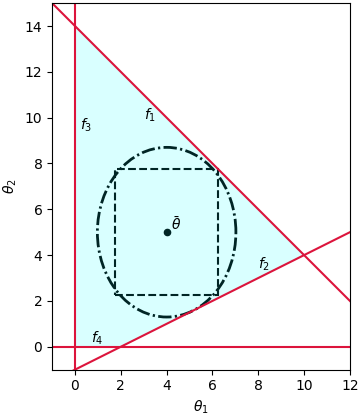}
		\caption{Case 2: $\textrm{Cov}(\theta_1,\theta_2) = 0$}
	\end{subfigure}
	\begin{subfigure}[t]{.3\textwidth}
		\centering
		\includegraphics[height=2.3in]{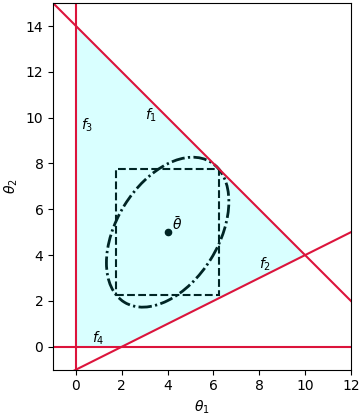}
		\caption{Case 3: $\textrm{Cov}(\theta_1,\theta_2) = 1$}
	\end{subfigure}
	\caption{A graphical depiction of system given by~\eqref{eq:2d_example_constraints} and~\eqref{eq:2d_example_conditions}. In each case the elliptical region (obtained with MICP approach) and traditional hyperbox region (obtained with MILP approach).}
	\label{fig:monte_carlo}
\end{figure}

Figure~\ref{fig:monte_carlo} shows the geometry of the cases detailed in Table~\ref{tab:2d_results} and includes the corresponding hyperbox uncertainty region dictated by \eqref{eq:hypercube} where $\Delta\boldsymbol{\theta}^-$ and $\Delta\boldsymbol{\theta}^+$ are both taken to be $(4.243, 5.196)$. These bounds correspond to $\bar{\theta}_i \pm 3\sigma_i$ where $\sigma_i$ is the standard deviation for the Gaussian random variable $\theta_i$. The flexibility index with hyperbox constraints was computed by solving the MILP formulation \eqref{eq:MILP_flexibility_index} and it was found to be $F=0.53$. In all three cases, it is visually apparent that inscribing an ellipsoidal region in the rectangular region will result in a smaller (suboptimal) ellipsoidal region relative to the one obtained via the proposed MICP method. Consequently, the MICP method is less conservative than the hyperbox approach. Also, it is interesting to note that the MICP approach identifies $f_2$ as the limiting constraint in cases 1 and 2, while the traditional hyperbox formulation always identifies $f_1$ as the limiting constraint.  This highlights the importance of capturing correlations.

\subsection{Heat Exchanger Network}
The heat exchanger network example presented in~\cite{grossmann1987active} is adapted for direct comparison with our approach. Figure~\ref{fig:hx_network} shows the network. The system constraints are:
\begin{equation}
	\begin{aligned}
		f_1 &= -350 K - 0.67 Q_c + T_3 \leq 0 \\
		f_2 &= 1388.5 K + 0.5 Q_c - 0.75 T_1 - T_3 - T_5 \leq 0 \\
		f_3 &= 2044 K + Q_c - 1.5 T_1 - 2 T_3 - T_5 \leq 0 \\
		f_4 &= 2830 K + Q_c - 1.5 T_1 - 2 T_3 - T_5 - 2 T_8 \leq 0 \\
		f_5 &= -3153 K - Q_c + 1.5 T_1 + 2 T_3 + T_5 + 3 T_8 \leq 0
	\end{aligned}
	\label{eq:hx_example_constraints}
\end{equation}
\noindent where $T_1$, $T_3$, $T_5$, and $T_8$ denote the Gaussian parameters and $Q_c$ denotes the recourse variable.  The mean and covariance matrix are given by:
\begin{equation}
	\bar{\boldsymbol{\theta}} =
	\begin{bmatrix}
		620 \\ 388 \\ 583 \\ 313
		\end{bmatrix} K
		\ \ \ \ \ \
		V_{\boldsymbol{\theta}} = 
		\begin{bmatrix}
		11.11 & \beta & \beta & \beta \\
		\beta & 11.11 & \beta & \beta \\
		\beta & \beta & 11.11 & \beta \\
		\beta & \beta & \beta & 11.11
	\end{bmatrix} K^2
	\label{eq:hx_example_conditions}
\end{equation}
\noindent {The parameter variance $\sigma_i^2$ is taken to be $11.11 K^2$ which corresponds to  $\bar{\theta}_i \pm 3 \sigma_i$, where $3 \sigma_i$ is equated to $10 K$ in accordance to the $\pm 10 K$ variations specified in \cite{grossmann1987active}.

\begin{figure}[ht]
	\centering
	\includegraphics[width=3.5in]{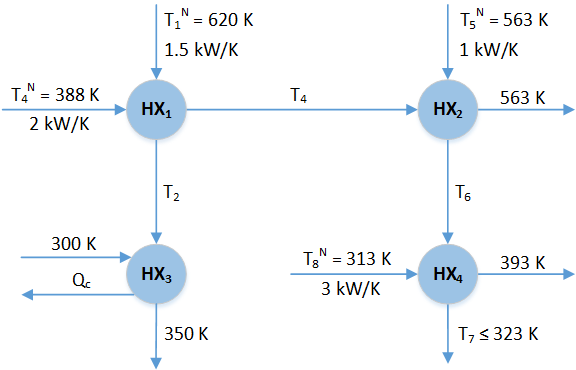}
	\caption{A heat exchanger network where $T_1$, $T_3$, $T_5$, and $T_8$ are Gaussian parameters and $Q_c$ is the recourse variable.}
	\label{fig:hx_network}
\end{figure}

We consider cases for $\beta=0$ and $\beta=5$ and the results are summarized in Table~\ref{tab:hx_results}. The stochastic flexibility index is again estimated via 100,000 MC samples of $\mathcal{N}(\bar{\boldsymbol{\theta}}, V_{\boldsymbol{\theta}})$ for each case. The $SF$-MC sampling requires $493$-$506$ seconds to estimate $SF$ while the MICP computation of $\alpha^*$ requires $0.15$-$0.18$ seconds. The {$\alpha^*$-MC estimation is again in agreement with the MICP results and we confirm that the confidence level is a lower bound for the stochastic flexibility index. We also note that the positive parameter correlation increases the system flexibility. Interestingly, in this case the gap between the confidence level and the stochastic flexibility index is quite large. This highlights that the shape of the feasible region and of the uncertainty set play a key role in the tightness of the lower bound provided by the proposed approach. 

\begin{table}[ht]
	\caption{Results for the 2 cases considered in connection to the system governed by~\eqref{eq:hx_example_constraints} and~\eqref{eq:hx_example_conditions}.}
	\begin{center}
		\begin{tabular}{l|ccccc}
			& $\beta$ & $\delta^*$ & $\alpha^*$-MICP (\%) & $\alpha^*$-MC (\%) &$SF$-MC (\%) \\ \hline
			Case 1 & 0 & 3.60 & 53.7 & 54.1 & 97.0\\
			Case 2 & 5 & 4.67 & 67.7 & 68.0 &97.1
		\end{tabular}
	\end{center}
	\label{tab:hx_results}
\end{table}

The hyperbox flexibility index given by problem~\eqref{eq:MILP_flexibility_index} is $F=0.5$. This means that the maximum allowable variations are $\pm 5 K$. This result can be juxtaposed against that of case 1 obtained with the MICP approach, which corresponds to a ellipsoidal uncertainty region with radius of $6.32 K$. Thus, the MICP formulation yields an ellipsoidal region that is appreciably larger (less conservative) than one that can be inscribed within the optimized hyperbox region. 

\begin{table}[!htp]
	\caption{Potential representations for uncertainty set $T(\delta)$.}
	\begin{center}
		\begin{tabular}{c|c}
			Name  & Uncertainty Set  \\ \hline
			Ellipsoidal-Norm & $T_{ellip}(\delta) = \{\boldsymbol{\theta} : ||\boldsymbol{\theta} - \bar{\boldsymbol{\theta}}||_{A}^2 \leq \delta \}$ \\
			$\ell_\infty$-Norm & $T_{\infty}(\delta) = \{\boldsymbol{\theta} : ||\boldsymbol{\theta} - \bar{\boldsymbol{\theta}}||_\infty \leq \delta \}$ \\
			$\ell_1$-Norm & $T_{1}(\delta) = \{\boldsymbol{\theta} : ||\boldsymbol{\theta} - \bar{\boldsymbol{\theta}}||_1 \leq \delta \}$ \\
			$\ell_2$-Norm & $T_{2}(\delta) = \{\boldsymbol{\theta} : ||\boldsymbol{\theta} - \bar{\boldsymbol{\theta}}||_2 \leq \delta \}$ \\
			CVaR-norm & $T_{CVaR}(\delta) = \{\boldsymbol{\theta} : \langle\langle\boldsymbol{\theta} - \bar{\boldsymbol{\theta}}\rangle\rangle_\alpha \leq \delta \}$
		\end{tabular}
	\end{center}
	\label{tab:bounded_norms}
\end{table}

\section{Conclusions and Future Work}

In this work we have presented a mixed-integer conic programming formulation to directly determine the flexibility index of a linear system in the face of multivariate Gaussian uncertainty. This approach is based on the active set method first proposed in \cite{grossmann1987active}. Here, we exploit the observation that any compact uncertainty set whose size can be measured in terms of a scalar value can be employed under this method. Our approach defines the flexibility index as the maximum radius of an  ellipsoidal uncertainty set and we show that the maximum radius can be used to compute the corresponding confidence level. We also showed that the confidence level provides a lower bound for the stochastic flexibility index (i.e., the probability of having feasible operation). The computational utility of this approach is currently being investigated in the context of large systems and will be explored further in future research. 

The observation that the flexibility index problem can use different uncertainty sets opens the possibility for a number of potential extensions. Table \ref{tab:bounded_norms} provides several examples of uncertainty set representations that can be used. A detailed analysis of the utility and mathematical properties of $\ell_p$-norm uncertainty sets in the context of robust optimization is provided in \cite{li2011comparative}. We note that these norms are scaled using the level set (as opposed to the linear scaling of hyperbox representations used in traditional flexibility analysis). This highlights that there exist alternatives to parameterize uncertainty sets. Moreover, we note that all these norms can be handled using existing modern mixed-integer programming solvers. In particular, the CVaR norm has recently received interest for its ability to approximate $\ell_p$-norms via linear programming formulations and because of this include the $\ell_1$ and $\ell_\infty$ norms as extreme cases \cite{gotoh2016two}. The utility of the different uncertainty sets in Table \ref{tab:bounded_norms} within the context of flexibility analysis has not yet been fully explored. 

\section*{Acknowledgments}

This work was supported by the U.S. Department of Energy under grant DE-SC0014114.

\appendix
\section{Proofs of Theorems}

\subsection{Proof of Theorem \ref{thm:active}}
\label{a:active}
Given that the set of gradients $\frac{\partial}{\partial \mathbf{z}} [f_j(\mathbf{z},\boldsymbol{\theta})], j \in J,$ are linearly independent, \eqref{eq:partial} requires that there exist $n > n_z$ nonzero Lagrange multipliers $\lambda_j$ since \eqref{eq:KKT_1} necessitates that at least one $\lambda_j$ be nonzero. The complete set of active constraints form a system of $n + n_\theta$ equations with $n_z + n_\theta + 1$ variables, thus a general solution requires $n \leq n_z + 1$. Hence, there must exist $n_z + 1$ active constraints at the solution to problem \eqref{eq:KKT_feasibility_function}. 

\subsection{Proof of Theorem \ref{thm:flexibility_index}}
\label{a:flexibility_index}
Suppose that there exists a bounded solution $\delta^*$, $\boldsymbol{\theta}^*$ to problem \eqref{eq:flexibility_index_definition} and that the system constraints $f_j(\mathbf{z}, \boldsymbol{\theta}), j \in J,$ are continuous in both $\mathbf{z}$ and $\boldsymbol{\theta}$. Now assume that the solution $\boldsymbol{\theta}^*$ is such that $\psi(\boldsymbol{\theta}^*) < 0$. For any $\delta^* < \hat{\delta}$ we have 

\begin{equation}
	\max_{\boldsymbol{\theta} \in T(\delta^*)} \psi(\boldsymbol{\theta}) = \psi(\boldsymbol{\theta}^*) \leq \psi(\boldsymbol{\hat{\theta}}) = \max_{\boldsymbol{\theta} \in T(\hat{\delta})} \psi(\boldsymbol{\theta}).
	\label{eq:proof1}
\end{equation}

\noindent Since $f_j(\mathbf{z}, \boldsymbol{\theta}), j \in J,$ are continuous in both $\mathbf{z}$ and $\boldsymbol{\theta}$, it follows that $\psi(\boldsymbol{\theta})$ is a continuous function in $\boldsymbol{\theta}$ (as shown in \cite{polak1979theoretical}). Hence, for some $\epsilon > 0$ and a $\tilde{\delta}$ in the neighborhood of the solution $\delta^*$ we have 

\begin{equation}
	|\psi(\boldsymbol{\theta}^*) - \psi(\boldsymbol{\tilde{\theta}})| < \epsilon,
	\label{eq:proof2}
\end{equation}
with $\psi(\boldsymbol{\tilde{\theta}}) = \max_{\boldsymbol{\theta} \in T(\tilde{\delta})} \psi(\boldsymbol{\theta})$. Rearranging \eqref{eq:proof2} we obtain

\begin{equation}
	\psi(\boldsymbol{\tilde{\theta}}) - \epsilon < \psi(\boldsymbol{\theta}^*) < 0
	\label{eq:proof3}
\end{equation}

\noindent and for $\epsilon$ sufficiently small it follows that

\begin{equation}
	\psi(\boldsymbol{\tilde{\theta}}) < 0.
	\label{eq:proof4}
\end{equation}

\noindent For a choice of $\tilde{\delta}$ arbitrary close to $\delta^*$ and which satisfies $\tilde{\delta} > \delta^*$, we have from \eqref{eq:proof1} and \eqref{eq:proof4} that $\tilde{\delta}$ is a feasible solution to problem \eqref{eq:flexibility_index_definition}. However, this is a contradiction to $\delta^*$ being a solution to \eqref{eq:flexibility_index_definition}. Thus, the assumption that $\psi(\boldsymbol{\theta}^*) < 0$ does not hold and $\psi(\boldsymbol{\theta}^*) = 0$.

Hence, for any compact set $T(\delta)$, the solution $\delta^*$ must lie on the boundary of the feasible region which is given by $\psi(\boldsymbol{\theta}) = 0$. Furthermore, the largest uncertainty set that can be fully inscribed in the feasible region is given by the smallest value of $\delta$ that satisfies $\psi(\boldsymbol{\theta}) = 0$. Thus problem \eqref{eq:flexibility_index_definition} can be reformulated as \eqref{eq:flexibility_index_mod}.

\bibliography{references}

\end{document}